\documentclass[a4paper,12pt]{article}
\usepackage{subcaption}
\usepackage{amsmath}
\usepackage{amssymb,eurosym}
\usepackage{mathtools}
\usepackage[dvips]{graphicx}
\usepackage{lscape}
\usepackage{color}
\usepackage{graphics}
\usepackage{graphicx}
\graphicspath{{./Graphiken/}}
\usepackage{amsmath}
\usepackage[nospace,noadjust]{cite}
\usepackage[pdftoolbar=true,
            pdfmenubar=true,
            pdfpagemode=UseOutlines,
            bookmarksnumbered=true,
            linktocpage=true,
            colorlinks=false]{hyperref}
\usepackage{verbatim}
\newtheorem{definition}{Definition}[section]
\newtheorem{theorem}[definition]{Theorem}
\newtheorem{lemma}[definition]{Lemma}
\newtheorem{corollary}[definition]{Corollary}
\newtheorem{proposition}[definition]{Proposition}
\newtheorem{remark}[definition]{Remark}
\newtheorem{example}[definition]{Example}

\newtheorem{assumption}[definition]{Assumption}

\newtheorem{conjecture}[definition]{Conjecture}

\newcommand{\bdefi}{\begin{definition}}
\newcommand{\edefi}{\end{definition}}
\DeclareMathOperator{\gem}{Gem}
\DeclareMathOperator{\new}{New}
\newcommand{\blem}{\begin{lemma}}
\newcommand{\elem}{\end{lemma}}
\newcommand{\bthe}{\begin{theorem}}
\newcommand{\ethe}{\end{theorem}}
\newcommand{\bcor}{\begin{corollary}}
\newcommand{\ecor}{\end{corollary}}
\newcommand{\bprop}{\begin{proposition}}
\newcommand{\eprop}{\end{proposition}}
\newcommand{\brem}{\begin{remark}}
\newcommand{\erem}{\end{remark}}
\newcommand{\bex}{\begin{example}}
\newcommand{\eex}{\end{example}}
\newcommand{\bass}{\begin{assumption}}
\newcommand{\eass}{\end{assumption}}

\newcommand{\bconj}{\begin{conjecture}}
\newcommand{\econj}{\end{conjecture}}

\renewcommand{\i}{\begin{itemize}}
\newcommand{\ii}{\end{itemize}}
\newcommand{\q}{\begin{equation}}
\newcommand{\qq}{\end{equation}}
\newcommand{\qa}{\begin{eqnarray*}}
\newcommand{\qqa}{\end{eqnarray*}}
\renewcommand{\a}{\begin{array}}
\renewcommand{\aa}{\end{array}}

\newcommand{\1}{{1\!\!1}}

\DeclareMathOperator{\conv}{conv}

\renewcommand{\epsilon}{\varepsilon}

\renewcommand{\L}{{\cal L}}

\newcommand{\N}{\mathbb N}

\renewcommand{\phi}{\varphi}
\newcommand{\pmat}[1]{\begin{pmatrix}#1\end{pmatrix}}

\newcommand{\proof}{{\bf Proof}. }
\newcommand{\qed}{~\hfill$\bullet$}
\newcommand{\R}{\mathbb R}

\renewcommand{\subset}{\subseteq}

\numberwithin{equation}{section}

\begin{document}

\title{\bf On Newton-polytope-type sufficiency conditions for coercivity of polynomials}

\author{
Tom\'a\v{s} Bajbar\thanks{Institute for Mathematics, Goethe University Frankfurt, Germany, tomas.bajbar@gmail.com}
\and
Yoshiyuki Sekiguchi\thanks{Institute for Mathematics, Tokyo University of Marine Science and Technology, Japan, yoshi-s@kaiyodai.ac.jp}
} 

%\date{preliminary version 1, \today}

\maketitle

\begin{abstract} We identify new sufficiency conditions for coercivity of general multivariate polynomials $f\in\mathbb{R}[x]$ which are expressed in terms of their Newton polytopes at infinity and which consist of a system of affine-linear inequalities in the space of polynomial coefficients. By sharpening the already existing necessary conditions for coercivity for a class of gem irregular polynomials we provide a characterization of coercivity of circuit polynomials, which extends the known results on this well studied class of polynomials. For the already existing sufficiency conditions for coercivity which contain a description involving a set projection operation, we identify an equivalent description involving a single posynomial inequality. This makes them more easy to apply and hence also more appealing from the practical perspective. We relate our results to the existing literature and we illustrate our results with several examples.
\end{abstract}



\newpage

\section{Introduction}\label{sec:intro} 

Let $\R[x] =\R[x_1,\dots,x_n]$ denote the ring of polynomials  in $n$ variables with real coefficients. A function $f:\R^n\rightarrow\R$ is called coercive on $\R^n$, if $f(x)\rightarrow +\infty$ holds whenever $\|x\|\rightarrow +\infty$, where $\| \cdot \|$ denotes some norm on $\R^n$. Coercivity of polynomials can be used as a sufficient condition for guaranteeing the existence of globally minimal points of optimization problems, which is often formulated as an assumption (see, e.g. \cite{Strm,Schwg,BS,AAA,Vui,Wang}). For analyzing the global invertibility property of polynomial maps $F:\R^n\to\R^n$ (see e.g. \cite{Biv,Tib}), an equivalence between the properness of $F$ and coercivity of $\|F\|_2^2\in\R[x]$ can be used for guaranteeing the global diffeomorphism property of $F$ (see e.g. \cite{BS2}). Since coercivity of $f\in\R[x]$ is equivalent to the boundedness of its lower level sets $\{x\in\R^n|\ f(x)\leq \alpha\}$ for all $\alpha\in\R$, appropriate coercivity sufficiency conditions are useful as a tool for analyzing the boundedness property of basic semialgebraic sets. 

This article is structured as follows. First the main concepts and results on coercive polynomials and their Newton polytopes at infinity from \cite{BS} are briefly summarized as they form the conceptual framework we shall work with in the present article.

In Section \ref{sec:coerc} with Theorem \ref{thm:main} we derive new sufficient conditions for coercivity of general polynomials which are expressed in terms of their Newton polytopes at infinity and which consist of a system of affine-linear inequalities in the space of polynomial coefficients. Then with Lemma \ref{lem:polyhed_suff} we prepare a reformulation (Theorem \ref{the:suffdeg_restated}) of the already existing coercivity sufficiency conditions (Theorem \ref{the:suffdeg}) which is easier to work with as it replaces the original projection-based formulation \eqref{eq:weigh1}-\eqref{eq:weigh3} by a single posynomial inequality \eqref{eq:the_suffdeg_restated}. With Theorem \ref{the:Char_Circ} we characterize coercivity of circuit polynomials and we also discuss the independence of coercivity sufficiency conditions from Theorems \ref{thm:main} and \ref{the:suffdeg_restated} in Examples \ref{ex:compare2} and \ref{ex:compare}. The article closes with some final remarks in Section \ref{sec:fin}.

In \cite{BS} it is shown how the coercivity of multivariate polynomials can often be analyzed by studying so-called Newton polytopes at infinity, whose definition we recall in the next step. 

We denote $\N_0:=\N\cup\{0\}$ and for $f\in\R[x]$ we write $f(x)=\sum_{\alpha\in A(f)}f_{\alpha}x^\alpha$ with $A(f)\subseteq\N_0^n$, $f_{\alpha}\in\R\setminus\{0\}$ for $\alpha\in A(f)$, and $x^{\alpha}=x_1^{\alpha_1}\dots x_n^{\alpha_n}$ for $\alpha\in\N_0^n$. For $f\in\R[x]$, the set $\new_{\infty}(f):=\conv \left(A(f)\cup\{0\}\right)$, that is, the convex hull of the set $A_0(f):=A(f)\cup\{0\}$ is called the \textit{Newton polytope at infinity} of the polynomial $f$, and the set $\new (f):=\conv \left(A(f)\right)$ is called the \textit{Newton polytope} of the polynomial $f$. The set $V_0(f)$ denotes the set of all vertices of $\new_{\infty}(f)$. With $\mathbb{H}:=\{h\in\R|\,h\geq 0\}$ one obtains, due to $0\in\new_{\infty}(f)\subseteq\mathbb{H}^n$, that the set $V_0(f)$ always contains the origin. For the later purposes of this work we shall define the vertex set of $\new_{\infty}(f)$ at infinity $V(f):=V_0(f)\setminus\{0\}$ and the set $V_0^c(f):=A(f)\setminus V_0(f)$ of all exponent vectors of $f$ which are no vertices of $\new_{\infty}(f)$ as well as the set $V^c(f):=A(f)\setminus V(f)$ of all exponent vectors of $f$ which are no vertices at infinity of $\new_{\infty}(f)$.

The following three conditions from \cite{BS} are crucial for analyzing the coercivity of $f\in\R[x]$ on $\R^n$. Here and subsequently we put $[n]:=\{1,\dots,n\}$.
\q
V(f)\subset 2\N_0^n.\tag{C1}\label{eq:c1}
\qq
\q
\text{All $\alpha\in V(f)$ satisfy $f_\alpha>0$.}\tag{C2}\label{eq:c2}
\qq
\q
\text{For all $i\in [n]$ the set $V(f)$ contains a vector of the form $2k_ie_i$ with $k_i\in\N$.}\tag{C3}\label{eq:c3}
\qq

For a polynomial $f\in\R[x]$ satisfying the condition \eqref{eq:c3}, we define $V_{\text{ess}}(f):=\{2k_ie_i,\ i\in [n]\}\subseteq V(f)$ to be the set of the \textit{essential} vertices at infinity of $\new_{\infty}(f)$.

For a polynomial $f\in\R[x]$ let
\[
\mathcal{G}(f)\ :=\ \{G\subset\R^n|\ G\neq\emptyset\ \text{is a face of}\ \new_{\infty}(f)\ \text{with}\ 0\not\in G\}
\]
be the set of all nonempty faces of  $\new_{\infty}(f)$ not including the origin. The set
\[
\gem(f)\ :=\ \bigcup_{G\in\mathcal{G}(f)}\,G,
\] 
is called the \emph{gem} of $f$ (in \cite{Kou, Pha08} also called the ``Newton boundary at infinity'') and gives rise to the following important regularity concept for polynomials:

\begin{definition}[\mbox{\cite{BS}}]\label{def:Fregular}  Let $f\in\R[x]$ be given.
\begin{itemize}
\item[a)] An exponent vector $\alpha\in A(f)$ is called \emph{gem degenerate} if  $\alpha\in V^c(f)\cap G$ holds for some $G\in\mathcal{G}(f)$. We denote the set of all gem degenerate points $\alpha\in A(f)$ by $D(f)$.
\item[b)] The polynomial $f$ is called \emph{gem re\-gular} if the set $D(f)$ is empty,
otherwise it is called \emph{gem irregular}. 
\end{itemize}
\end{definition}

We recall the following characterization of the set $D(f)$ from \cite{BS} which states that $D(f)$ contains exactly the exponent vectors in $A(f)$ which cannot be written as convex combination of elements from $V_0(f)$ with the origin entering with a positive weight.

\bprop[Characterization of the set $D(f)$, \mbox{\cite[Prop. 2.24]{BS}}]\label{prop:BS} For a polynomial $f\in\R[x]$ satisfying the condition \eqref{eq:c3} the following are equivalent:
\begin{itemize}
\item[(a)] $\alpha^\star\in D(f)$
\item[(b)] $\alpha^\star\in V^c(f)$, and any choice of coefficients $\lambda_\alpha$, $\alpha\in V_0(f)$ with
\[
\alpha^\star=\sum_{\alpha\in V_0(f)}\lambda_\alpha \alpha,\quad\sum_{\alpha\in V_0(f)}\lambda_\alpha=1,\quad\lambda_\alpha\geq 0,\ \alpha\in V_0(f)
\]
satisfies $\lambda_0=0.$
\end{itemize}
\eprop

Clearly, gem regularity of $f\in\R[x]$ is equivalent to $V^c(f)\cap G=\emptyset$ for all $G\in\cal{G}$.
Furthermore, the definition of $D(f)$ gives rise to a partitioning of $V^c(f)$ into $D(f)$ and a set of ``remaining exponents''
$R(f):=V^c(f)\setminus D(f)$, so that we may write
\q\label{eq:VDR}
A(f)\ =\ V(f)\ \dot\cup\ D(f)\ \dot\cup\ R(f).
\qq
Using \eqref{eq:VDR} together with the notation $f^S(x):=\sum_{\alpha\in S}f_\alpha x^\alpha$ for some $S\subseteq A(f)$, any $f\in\R[x]$ can be expressed as
\q\label{eq:fVDR}
f=f^{V(f)}+f^{D(f)}+f^{R(f)}.
\qq

Now we are ready to state the general necessary conditions for coercivity of polynomials.

\bthe[Necessary conditions for coercivity \mbox{\cite[Th. 2.8]{BS}}]\label{the:necessary_general}\ \\ Let $f\in\R[x]$ be coercive on $\R^n$. Then $f$ fulfills the conditions \eqref{eq:c1}--\eqref{eq:c3}.
\ethe

Although the coniditions \eqref{eq:c1}--\eqref{eq:c3} are not sufficient for coercivity of polynomials in general, in the following theorem it was shown that they are sufficient for coercivity for a broad class of gem regular polynomials.

\bthe[Characterization of coercivity, \mbox{\cite[Th. 3.2]{BS}}]\label{the:char}\ \\
Let $f\in\R[x]$ be gem regular. Then the following assertions are equivalent.

\i
\item[a)] $f$ is coercive on $\R^n$.
\item[b)] $f$ fulfills the conditions \eqref{eq:c1}-\eqref{eq:c3}.
\ii
\ethe

Recall that, by Carath\'eodory's theorem, for any exponent vector $\alpha^\star\in V_0^c(f)$ there exists a set of affinely independent points
$W_{\alpha^\star}\subset V_0(f)$ with $\alpha^\star\in\conv W_{\alpha^\star}$. In the case that a simplicial face $G$ of $\new_{\infty}(f)$ contains $\alpha^\star$, the
set $W_{\alpha^\star}$ can be chosen as the vertex set $\text{Vert}(G)$ of $G$. For non-simplicial faces $G$, however, there may exist several possibilities
to choose $W_{\alpha^\star}\subset \text{Vert}(G)$. If, in addition, $W_{\alpha^\star}$ is chosen minimally in the sense that the presence of all points in 
$W_{\alpha^\star}$ is necessary for $\alpha^\star\in\conv W_{\alpha^\star}$ to hold, then 
we also have $\lambda_\alpha>0$ for all $\alpha\in W_{\alpha^\star}$. Note that a minimal choice of $W_{\alpha^\star}$ is not necessarily unique. This idea of 'minimality' gives rise to the following definition which proves to be convenient for the purposes of the present article.

\begin{definition}[Map of minimal barycentric coordinates]\label{def:bary_map} For a given polynomial $f\in\R[x]$ a map $\lambda_f:V_0^c(f)\times V_0(f)\to [0,1]$ is called a map of minimal barycentric coordinates of $f$, if for each $\alpha^\star\in V_0^c(f)$ there exists some affinely independent set $W_{\alpha^\star}\subseteq V_0(f)$ such that
\begin{itemize}
\item[i)]$\lambda_f(\alpha^\star,\alpha)>0$ for all $\alpha\in W_{\alpha^\star}$
\item[ii)] $\lambda_f(\alpha^\star,\alpha)=0$ for all $\alpha\in V_0(f)\setminus W_{\alpha^\star}$
\item[iii)] $
\sum_{\alpha\in W_{\alpha^\star}}\lambda_f(\alpha^\star,\alpha)\pmat{\alpha\\1}\ =\ \pmat{\alpha^\star\\1}
$
\end{itemize}
hold. Given a map of minimal barycentric coordinates $\lambda_f$ of some polynomial $f\in\R[x]$ and given some $\alpha^\star\in V_0^c(f)$, we call the set of affinely independent points $W_{\alpha^\star}(\lambda_f):=\{\alpha\in V_0(f)|\ \lambda_f(\alpha^\star,\alpha)>0\}$ the minimal vertex representation of $\alpha^\star$ corresponding to $\lambda_f$.
\end{definition}

Since for any set of affinely independent points $W_{\alpha^\star}$ with $\alpha^\star\in\conv W_{\alpha^\star}$, the set of solutions $\{\lambda_f(\alpha^\star,\alpha)\in\R,\ \alpha\in W_{\alpha^\star}\}$ corresponding to the system
\[%\label{eq:lambdaVstar}
\sum_{\alpha\in W_{\alpha^\star}}\lambda_f(\alpha^\star,\alpha)\pmat{\alpha\\1}\ =\ \pmat{\alpha^\star\\1},\quad \lambda_f(\alpha^\star,\alpha)\ge0,\ \alpha\in W_{\alpha^\star}
\]
is unique, then for any \textit{fixed} choice of the subsets $W_{\alpha^\star}$ in Definition \ref{def:bary_map}, also the corresponding map of minimal barycentric coordinates $\lambda_f$ is uniquely determined. Clearly, for a general polynomial $f\in\R[x]$ there might exist several maps of minimal barycentric coordinates $\lambda_f$.

\brem\label{rem:Bary_Deg} In view of Proposition \ref{prop:BS}, for any polynomial $f\in\R[x]$ satisfying the condition \eqref{eq:c3} and any map of minimal barycentric coordinates $\lambda_f$ of $f$, one obtains that $0\notin W_{\alpha^\star}(\lambda_f)$, and consequently, $W_{\alpha^\star}(\lambda_f)\subseteq V(f)$ holding for all $\alpha^\star\in D(f)$.
\erem

For any polynomial $f\in\R[x]$ and any map of minimal barycentric coordinates $\lambda_f$ of $f$, we may consider for each $\alpha^\star\in V_0^c(f)$ the \emph{circuit number} (cf. \cite{Ili14Wa}) 
\[
\Theta(f,\lambda_f,\alpha^\star)\ =\ \prod_{\alpha\in W_{\alpha^\star}(\lambda_f)}\left(\frac{f_\alpha}{\lambda_f(\alpha^\star,\alpha)}\right)^{\lambda_f(\alpha^\star,\alpha)}.
\]

\brem\label{rem:circ_positive}
Clearly, for any polynomial $f\in\R[x]$ satisfying the conditions \eqref{eq:c2} and \eqref{eq:c3} and any map of minimal barycentric coordinates $\lambda_f$ of $f$, the circuit number $\Theta(f,\lambda_f,\alpha^\star)$ corresponding to each gem degenerate exponent vector $\alpha^\star\in D(f)$ is positive: In fact, by \eqref{eq:c3} and Remark \ref{rem:Bary_Deg} one has $0\notin W_{\alpha^\star}(\lambda_f)\subseteq V(f)$ and condition \eqref{eq:c2} implies $f_\alpha>0$ for each $\alpha\in W_{\alpha^\star}(\lambda_f)$. 
\erem

We will see that for guaranteeing coercivity of some polynomial $f$ only those circuit numbers corresponding to the gem degenerate exponent vectors $\alpha^\star\in D(f)$ are important. We recall the following result from \cite{BS} which, unlike Theorem~\ref{the:char}, guarantees coercivity even for a broad class of gem irregular polynomials. We restate this result here by using the map of minimal barycentric coordinates from Definition \ref{def:bary_map}.

\bthe[Sufficient conditions for coercivity, \mbox{\cite[Th. 3.4]{BS}}]\label{the:suffdeg}\ \\
Let $f\in\R[x]$ be a polynomial satisfying the conditions \eqref{eq:c1}-\eqref{eq:c3} and let $\lambda_f$ be some map of minimal barycentric coordinates of $f$. Furthermore, for each $\alpha^\star\in D(f)$ let $\omega(\alpha^\star)>0$ denote weights such that 
\[\sum_{\alpha^\star\in D(f)}\omega(\alpha^\star)\le1\] 
holds and let
\[
f_{\alpha^\star}\ >\ -\omega(\alpha^\star)\,\Theta(f,\lambda_f,\alpha^\star)\ \text{ if }\ \alpha^\star\in 2\N_0^n
\]
and 
\[
|f_{\alpha^\star}|\ <\ \omega(\alpha^\star)\,\Theta(f,\lambda_f,\alpha^\star)\ \text{ else.}
\] 

Then $f$ is coercive on $\R^n$.
\ethe

Finally we recall the definition of a circuit polynomial from \cite{Ili14Wa,DIdW}. 

\bdefi[Circuit polynomial \mbox{\cite[Def. 2.2]{DIdW}}]\label{defi:circuit_polynomial} A polynomial $f\in\R[x]$ of the form 
\[
f(x)=\sum_{j=0}^rf_{\alpha(j)}x^{\alpha(j)}+f_{\alpha^\star}x^{\alpha^\star}
\]
with $r\leq n$ is called a circuit polynomial if the following conditions are fulfilled:
\begin{itemize}
\item[i)] $\alpha(j)\in 2\N_0^n$ for all $j=0,\dots,r$
\item[ii)] $f_{\alpha(j)}>0$ for all $j=0,\dots,r$
\item[iii)] $\text{Vert}\,(\new(f))=\{\alpha(0),\alpha(1),\dots,\alpha(r)\}$ with $\alpha(0),\alpha(1),\dots,\alpha(r)$ affinely independent.
\item[iv)] the exponent $\alpha^\star$ can be written uniquely as
\[
\alpha^\star=\sum_{j=0}^r\lambda_j\alpha(j)\quad\text{with}\quad \lambda_j>0\quad\text{and}\quad \sum_{j=0}^r\lambda_j=1
\]
in barycentric coordinates $\lambda_j$ relative to the vertices $\alpha(j)$, $j=0,\dots,r$.
\end{itemize} 
\edefi

For circuit polynomials the following characterization of their global non-negativity via circuit numbers was shown (see, e.g. \cite{Ili14Wa,DIdW}) which was recently successfully used within the polynomial optimization area for developing the Sum-of-Nonnegative-Circuits-based (abbr. SONC) algorithmic solution approach (for more details see, e.g. \cite{DIdW, Wang2, Murray, Theobald}).

\bthe[Non-negativity of circuit polynomials \mbox{\cite[Th. 2.3]{DIdW}}]\label{thm:IliWolff} Let $f\in\R[x]$ be a circuit polynomial. Then $f(x)\geq 0$ for all $x\in\R^n$ if and only if
\begin{equation}\label{eq:Circ1}
f_{\alpha^\star}\geq -\prod_{j=0}^r\left(\frac{f_{\alpha(j)}}{\lambda_j}\right)^{\lambda_j}\quad\text{ if }\alpha^\star\in 2\N_0^n
\end{equation}
and
\begin{equation}\label{eq:Circ2}
|f_{\alpha^\star}|\leq \prod_{j=0}^r\left(\frac{f_{\alpha(j)}}{\lambda_j}\right)^{\lambda_j}\quad\text{ if }\alpha^\star\notin 2\N_0^n.
\end{equation}
\ethe

\section{Main results}\label{sec:coerc}

The crucial part in the proof of the Characterization Theorem \ref{the:char} turns out to be the result asserting that the growth of gem regular polynomials $f\in\R[x]$ at infinity is governed from below by the part of the polynomial $f$ corresponding to its vertices at infinity $V(f)$. We restate this result briefly as we will use it for the proof of our new coercivity sufficiency conditions in Theorem \ref{thm:main}.

\blem[\mbox{\cite[Prop. 3.1]{BS}}]
\label{lem:basic_ineq}
Let $f$ be a gem regular polynomial satisfying the conditions
\eqref{eq:c1}--\eqref{eq:c3}. Then 
for any sequence of points $(x^k)_{k\in\N}$ with $\lim_{k\to +\infty}\|x^k\|=+\infty$ there exists some $\epsilon>0$ with 
\begin{equation}
f(x^k)\geq \epsilon f^{V(f)}(x^k) \quad\text{for almost all }k\in\N.
\end{equation}
\elem

The following theorem identifies polyhedral subsets in the space of coefficients which guarantee the coercivity property for a broad class of polynomials. 

\bthe [Sufficient condition for coercivity]\label{thm:main} Let $f\in\R[x]$ be a polynomial satisfying the conditions \eqref{eq:c1}-\eqref{eq:c3} and let $\lambda_f$ be a map of minimal barycentric coordinates of $f$. If the inequality

\begin{equation}\label{eq:main}
f_{\alpha}>\sum_{\alpha^\star\in D(f)\,\cap\, (2\N_0^n)^c}|f_{\alpha^\star}|\,\lambda_f(\alpha^\star,\alpha) -\sum_{\alpha^\star\in D(f)\,\cap\, 2\N_0^n}\min\{0,f_{\alpha^\star}\}\,\lambda_f(\alpha^\star,\alpha)
\end{equation}
is satisfied for each $\alpha\in V(f)$, then $f$ is coercive on $\R^n$.
\ethe

\proof According to \eqref{eq:VDR} and \eqref{eq:fVDR} we can write
$f=f^{V(f)}+f^{D(f)}+f^{R(f)}$.

First we define
\begin{equation}\label{eq:h_alpha}
h_{\alpha}:=\sum_{\alpha^\star\in D(f)\,\cap\, (2\N_0^n)^c}|f_{\alpha^\star}|\lambda_f(\alpha^\star,\alpha)\quad -\sum_{\alpha^\star\in D(f)\,\cap\, 2\N_0^n}\text{min}\{0,f_{\alpha^\star}\}\lambda_f(\alpha^\star,\alpha)
\end{equation}
for each $\alpha\in V(f)$. Notice that one can write
\[
f^{V(f)}(x)+f^{D(f)}(x)=F(x)+G(x)
\]
with
\[
F(x):=\sum_{\alpha\in V(f)}\left(f_{\alpha}-h_{\alpha}\right)x^\alpha\qquad \text{ and }\qquad G(x):=f^{D(f)}(x)+\sum_{\alpha\in V(f)}h_\alpha x^\alpha.
\]

If we show that $F$ is coercive and $G$ is globally
non-negative on $\R^n$, then we have coercivity of $f$.
In fact, since $F + f^{R(f)}$ is gem regular
and satisfies \eqref{eq:c1}--\eqref{eq:c3},
Lemma \ref{lem:basic_ineq} implies that
for any sequence of points $(x^k)_{k\in\N}$ with $\lim_{k\to +\infty}\|x^k\|=+\infty$ there exists some $\epsilon>0$ with 
\begin{equation}
 F(x^k) + f^{R(f)}(x^k)\geq \epsilon F(x^k) \quad\text{for almost all
  }k\in\N.
  \label{eq:suff_1}
\end{equation}
Since global non-negativity of $G$ and \eqref{eq:suff_1} imply
\begin{align*}
 f(x^k) & = f^{V(f)}(x^k) + f^{D(f)}(x^k) + f^{R(f)}(x^k) \\
 & = F(x^k) + G(x^k) + f^{R(f)}(x^k) \\
 & \geq F(x^k) + f^{R(f)}(x^k) \geq \epsilon F(x^k).
\end{align*} 
holding for almost all $k\in \N$, 
the coercivity of $F$ on $\R^n$ yields $\lim_{k\to+\infty}f(x^k)=+\infty$, which concludes the proof.

In order to show the coercivity of $F$ on $\R^n$, we partition the set of vertices at infinity $V(f)$ of $f$ into the set of the so-called essential vertices $V_{\text{ess}}(f)\subseteq V(f)$ corresponding to vertices at infinity of $f$ which lie on the axes of $\R^n$ as defined in the first chapter, and, into the remaining vertices at infinty of $f$ in the following sense:
\[
V(f)=V_{\text{ess}}(f)\ \dot\cup\ (V(f)\setminus V_{\text{ess}}(f)).
\]
This yields the decomposition
\begin{equation}\label{eq:decompo}
F(x)=\sum_{\alpha\in V_{\text{ess}}(f)}\left(f_{\alpha}-h_{\alpha}\right)x^\alpha\quad +\sum_{\alpha\in V(f)\setminus V_{\text{ess}}(f)}\left(f_{\alpha}-h_{\alpha}\right)x^\alpha.
\end{equation}
Notice first, that due to the assumption \eqref{eq:main} one has $f_\alpha-h_\alpha>0$ for all $\alpha\in V(f)$. Since $f$ satisfies the conditions \eqref{eq:c1}-\eqref{eq:c3}, the latter implies that the polynomial
\[
\sum_{\alpha\in V_{\text{ess}}(f)}\left(f_{\alpha}-h_{\alpha}\right)x^\alpha
\]
is coercive on $\R^n$ and is further implies that the polynomial
\[
\sum_{\alpha\in V(f)\setminus V_{\text{ess}}(f)}\left(f_{\alpha}-h_{\alpha}\right)x^\alpha
\]
is globally non-negative on $\R^n$. Then by \eqref{eq:decompo} the polynomial $F$ is coercive on $\R^n$ as a sum of a coercive polynomial and a globally non-negative one.

In order to prove the global non-negativity of $G$ on $\R^n$ we first show that for each $\alpha^\star\in D(f)$ one has
\begin{equation}\label{eq:g_nonneg}
g_{\alpha^\star}(x):=\sum_{\alpha\in V_0(f)}|f_{\alpha^\star}|\,\lambda_f(\alpha^\star,\alpha)\,x^\alpha+f_{\alpha^\star}x^{\alpha^\star}\geq 0\quad\text{for all }x\in\R^n.
\end{equation}
Let $\alpha^\star\in D(f)$. Due to Definition \ref{def:bary_map} of $\lambda_f$ and Remark \ref{rem:Bary_Deg}, there exists some subset $W_{\alpha^\star}(\lambda_f)\subseteq V(f)$ of affinely independent vertices at infinity of $\new_{\infty}(f)$ such that $\lambda_f(\alpha^\star,\alpha)>0$ for all $\alpha\in W_{\alpha^\star}(\lambda_f)$, $\lambda_f(\alpha^\star,\alpha)=0$ for all $\alpha\in V_0(f)\setminus W_{\alpha^\star}(\lambda_f)$ and $\alpha^\star=\sum_{\alpha\in W_{\alpha^\star}(\lambda_f)}\lambda_f(\alpha^\star,\alpha)\,\alpha$ with $\sum_{\alpha\in W_{\alpha^\star}(\lambda_f)}\lambda_f(\alpha^\star,\alpha)=1$. This yields
\[
g_{\alpha^\star}(x)=\sum_{\alpha\in W_{\alpha^\star}(\lambda_f)}|f_{\alpha^\star}|\,\lambda_f(\alpha^\star,\alpha)\,x^\alpha+f_{\alpha^\star}x^{\alpha^\star}
\]
and, since $W_{\alpha^\star}(\lambda_f)\subseteq 2\N_0^n$ holds due to condition \eqref{eq:c2}, $g_{\alpha^\star}$ is a circuit polynomial as defined in Theorem \ref{thm:IliWolff}. For the circuit polynomial $g_{\alpha^\star}$ one obtains
\[
\prod_{\alpha\in W_{\alpha^\star}(\lambda_f)}\left(\frac{|f_{\alpha^\star}|\lambda_f(\alpha^\star,\alpha)}{\lambda_f(\alpha^\star,\alpha)}\right)^{\lambda_f(\alpha^\star,\alpha)}=\prod_{\alpha\in W_{\alpha^\star}(\lambda_f)}|f_{\alpha^\star}|^{\lambda_f(\alpha^\star,\alpha)}
\]
\[
=|f_{\alpha^\star}|^{\sum_{\alpha\in W_{\alpha^\star}(\lambda_f)}\lambda_f(\alpha^\star,\alpha)}=|f_{\alpha^\star}|,
\]
and a direct application of Theorem \ref{thm:IliWolff} by checking the conditions \eqref{eq:Circ1}-\eqref{eq:Circ2} reveals that $g_{\alpha^\star}(x)\geq 0$ for all $x\in\R^n$, which proves \eqref{eq:g_nonneg}.

Finally, for the polynomial $G$ the following non-negativity estimate holds
\[
G(x)=f^{D(f)}(x)+\sum_{\alpha\in V(f)}h_\alpha x^\alpha
\]
\[
= \sum_{\alpha\in V(f)}\left( \sum_{\alpha^\star\in D(f)\,\cap\, (2\N_0^n)^c}|f_{\alpha^\star}|\lambda_f(\alpha^\star,\alpha)\, -\sum_{\alpha^\star\in D(f)\,\cap\, 2\N_0^n}\text{min}\{0,f_{\alpha^\star}\}\lambda_f(\alpha^\star,\alpha) \right) x^\alpha
\]
\[ + \sum_{\alpha^\star\in D(f)\,\cap\, (2\N_0^n)^c}f_{\alpha^\star}x^{\alpha^\star} + \sum_{\alpha^\star\in D(f)\, \cap\, 2\N_0^n}f_{\alpha^\star}x^{\alpha^\star}
\]
\[
\geq \sum_{\alpha\in V(f)}\left( \sum_{\alpha^\star\in D(f)\,\cap\, (2\N_0^n)^c}|f_{\alpha^\star}|\lambda_f(\alpha^\star,\alpha)\, -\sum_{\alpha^\star\in D(f)\,\cap\, 2\N_0^n}\text{min}\{0,f_{\alpha^\star}\}\lambda_f(\alpha^\star,\alpha) \right) x^\alpha
\]
\[ + \sum_{\alpha^\star\in D(f)\,\cap\, (2\N_0^n)^c}f_{\alpha^\star}x^{\alpha^\star} + \sum_{\substack{\alpha^\star\in D(f)\, \cap\, 2\N_0^n \\ f_{\alpha^\star} <0}}f_{\alpha^\star}x^{\alpha^\star}
\]
\[
= \sum_{\alpha\in V(f)}\left( \sum_{\alpha^\star\in D(f)\,\cap\, (2\N_0^n)^c}|f_{\alpha^\star}|\lambda_f(\alpha^\star,\alpha)\, -\sum_{\substack{\alpha^\star\in D(f)\, \cap\, 2\N_0^n \\ f_{\alpha^\star}<0}}\text{min}\{0,f_{\alpha^\star}\}\lambda_f(\alpha^\star,\alpha) \right) x^\alpha
\]
\[ + \sum_{\alpha^\star\in D(f)\,\cap\, (2\N_0^n)^c}f_{\alpha^\star}x^{\alpha^\star} + \sum_{\substack{\alpha^\star\in D(f)\, \cap\, 2\N_0^n \\ f_{\alpha^\star} <0}}f_{\alpha^\star}x^{\alpha^\star}
\]

\[
= \sum_{\alpha\in V(f)}\ \sum_{\alpha^\star\in D(f)\,\cap\, (2\N_0^n)^c} |f_{\alpha^\star}|\lambda_f(\alpha^\star,\alpha) \, x^\alpha \quad + \sum_{\alpha^\star\in D(f)\,\cap\, (2\N_0^n)^c}f_{\alpha^\star}x^{\alpha^\star}
\]
\[ 
+\sum_{\alpha\in V(f)}\  \sum_{\substack{\alpha^\star\in D(f)\, \cap\, 2\N_0^n \\ f_{\alpha^\star}<0}}-\text{min}\{0,f_{\alpha^\star}\}\lambda_f(\alpha^\star,\alpha)\,  x^\alpha\quad + \sum_{\substack{\alpha^\star\in D(f)\, \cap\, 2\N_0^n \\ f_{\alpha^\star} <0}}f_{\alpha^\star}x^{\alpha^\star}
\]

\[
= \sum_{\alpha^\star\in D(f)\,\cap\, (2\N_0^n)^c}\  \sum_{\alpha\in V(f)}|f_{\alpha^\star}|\lambda_f(\alpha^\star,\alpha) \, x^\alpha \quad + \sum_{\alpha^\star\in D(f)\,\cap\, (2\N_0^n)^c}f_{\alpha^\star}x^{\alpha^\star}
\]
\[ 
+\sum_{\substack{\alpha^\star\in D(f)\, \cap\, 2\N_0^n \\ f_{\alpha^\star}<0}}\ \sum_{\alpha\in V(f)}\, -\text{min}\{0,f_{\alpha^\star}\}\lambda_f(\alpha^\star,\alpha)\,  x^\alpha\quad + \sum_{\substack{\alpha^\star\in D(f)\, \cap\, 2\N_0^n \\ f_{\alpha^\star} <0}}f_{\alpha^\star}x^{\alpha^\star}
\]

\[
= \sum_{\alpha^\star\in D(f)\,\cap\, (2\N_0^n)^c}\left(  \sum_{\alpha\in V(f)}|f_{\alpha^\star}|\lambda_f(\alpha^\star,\alpha) \, x^\alpha \quad + f_{\alpha^\star}x^{\alpha^\star}\right)
\]
\[ 
+\sum_{\substack{\alpha^\star\in D(f)\, \cap\, 2\N_0^n \\ f_{\alpha^\star}<0}}\left( \sum_{\alpha\in V(f)}\, -\text{min}\{0,f_{\alpha^\star}\}\lambda_f(\alpha^\star,\alpha)\,  x^\alpha\quad + f_{\alpha^\star}x^{\alpha^\star}\right)
\]
\[
= \sum_{\alpha^\star\in D(f)\,\cap\, (2\N_0^n)^c}\left(  \sum_{\alpha\in V(f)}|f_{\alpha^\star}|\lambda_f(\alpha^\star,\alpha) \, x^\alpha \quad + f_{\alpha^\star}x^{\alpha^\star}\right)
\]
\[ 
+\sum_{\substack{\alpha^\star\in D(f)\, \cap\, 2\N_0^n \\ f_{\alpha^\star}<0}}\left( \sum_{\alpha\in V(f)}\, |f_{\alpha^\star}|\lambda_f(\alpha^\star,\alpha)\,  x^\alpha\quad + f_{\alpha^\star}x^{\alpha^\star}\right)
\]
\[
= \sum_{\alpha^\star\in D(f)\,\cap\, (2\N_0^n)^c} g_{\alpha^\star}(x)+ \sum_{\substack{\alpha^\star\in D(f)\, \cap\, 2\N_0^n \\ f_{\alpha^\star}<0}} g_{\alpha^\star}(x)\qquad  \geq\qquad  0.
\]
Thus the polynomial $G$ is globally non-negative on $\R^n$.
\qed

Next we turn our attention to the context of Theorem \ref{the:suffdeg}. The following lemma offers a way to reformulate Theorem \ref{the:suffdeg} by eliminating the appearing weights $\omega(\alpha^\star),\ \alpha^\star\in D(f)$ by using a single posynomial inequality. 

\blem\label{lem:polyhed_suff} Let all the assumptions and notation from Theorem \ref{the:suffdeg} be given. Then there exist weights $\omega(\alpha^\star)>0$, $\alpha^\star\in D(f)$ satisfying 
\begin{equation}\label{eq:weigh1}
\sum_{\alpha^\star\in D(f)}\omega(\alpha^\star)\le1\end{equation} 
\begin{equation}\label{eq:weigh2}
f_{\alpha^\star}\ >\ -\omega(\alpha^\star)\,\Theta(f,\lambda_f,\alpha^\star)\ \text{ if }\ \alpha^\star\in 2\N_0^n
\end{equation} 
\begin{equation}\label{eq:weigh3}
|f_{\alpha^\star}|\ <\ \omega(\alpha^\star)\,\Theta(f,\lambda_f,\alpha^\star)\ \text{ if }\ \alpha^\star\in (2\N_0^n)^c
\end{equation}
if and only if
\begin{equation}\label{eq:polyhed_suff}
\sum_{\alpha^\star\in D(f)\,\cap\, (2\N_0^n)^c}\frac{|f_{\alpha^\star}|}{\Theta(f,\lambda_f,\alpha^\star)} -\sum_{\alpha^\star\in D(f)\, \cap\, 2\N_0^n }\frac{\min\{0,f_{\alpha^\star}\}}{\Theta(f,\lambda_f,\alpha^\star)}<1
\end{equation}
\elem

\proof
Under the assumptions of Theorem \ref{the:suffdeg}, the polynomial $f$ satisfies the conditions \eqref{eq:c1}-\eqref{eq:c3} and according to the Remark \ref{rem:circ_positive} the expression \eqref{eq:polyhed_suff} is well defined since all appearing circuit numbers $\Theta(f,\lambda_f,\alpha^\star)$ are positive. 

For the first direction let some weights $\omega(\alpha^\star)>0$, $\alpha^\star\in D(f)$ be given which satisfy \eqref{eq:weigh1}-\eqref{eq:weigh3}. Then by \eqref{eq:weigh3}
\begin{equation}\label{eq:weigh_proof1}
\sum_{\alpha^\star\in D(f)\,\cap\, (2\N_0^n)^c}\frac{|f_{\alpha^\star}|}{\Theta(f,\lambda_f,\alpha^\star)}<\sum_{\alpha^\star\in D(f)\,\cap\, (2\N_0^n)^c}\omega(\alpha^\star).
\end{equation}
Further, due to $\min\{0,f_{\alpha^\star}\}=0$ for $f_{\alpha^\star}\geq 0$, -$\min\{0,f_{\alpha^\star}\}=-f_{\alpha^\star}$ for $f_{\alpha^\star}<0$ and the property \eqref{eq:weigh2} one obtains
\[
\sum_{\alpha^\star\in D(f)\, \cap\, 2\N_0^n }\frac{-\min\{0,f_{\alpha^\star}\}}{\Theta(f,\lambda_f,\alpha^\star)}
\]
\begin{equation}\label{eq:weigh_proof2}
= \sum_{\substack{\alpha^\star\in D(f)\, \cap\, 2\N_0^n \\ f_{\alpha^\star} <0}}\frac{-f_{\alpha^\star}}{\Theta(f,\lambda_f,\alpha^\star)}<\sum_{\substack{\alpha^\star\in D(f)\, \cap\, 2\N_0^n \\ f_{\alpha^\star} <0}}\omega(\alpha^\star).
\end{equation}
Combining \eqref{eq:weigh_proof1} and \eqref{eq:weigh_proof2} and applying \eqref{eq:weigh1} yields
\[
\sum_{\alpha^\star\in D(f)\,\cap\, (2\N_0^n)^c}\frac{|f_{\alpha^\star}|}{\Theta(f,\lambda_f,\alpha^\star)}+\sum_{\alpha^\star\in D(f)\, \cap\, 2\N_0^n }\frac{-\min\{0,f_{\alpha^\star}\}}{\Theta(f,\lambda_f,\alpha^\star)} 
\]
\[
< \sum_{\alpha^\star\in D(f)\,\cap\, (2\N_0^n)^c}\omega(\alpha^\star) + \sum_{\substack{\alpha^\star\in D(f)\, \cap\, 2\N_0^n \\ f_{\alpha^\star} <0}}\omega(\alpha^\star)\leq \sum_{\alpha^\star\in D(f)}\omega(\alpha^\star)\quad \leq\quad 1
\]
which implies \eqref{eq:polyhed_suff}. 

For the other direction let \eqref{eq:polyhed_suff} be satisfied. Then there exists some $\epsilon>0$ such that 
\begin{equation}\label{eq:epsilon}
\sum_{\alpha^\star\in D(f)\,\cap\, (2\N_0^n)^c}\frac{|f_{\alpha^\star}|}{\Theta(f,\lambda_f,\alpha^\star)} -\sum_{\alpha^\star\in D(f)\, \cap\, 2\N_0^n }\frac{\min\{0,f_{\alpha^\star}\}}{\Theta(f,\lambda_f,\alpha^\star)}=1-\epsilon
\end{equation} In order prove the existence of some weights $\omega(\alpha^\star)$, $\alpha^\star\in D(f)$ which satisfy \eqref{eq:weigh1}-\eqref{eq:weigh3} we first introduce the sets
\[
\Delta_1:=\{\alpha^\star\in D(f):\ \alpha^\star\in (2\N_0^n)^c \text{ or } \{\alpha^\star\in (2\N_0^n) \text{ with }f_{\alpha^\star}<0\}\}
\]
and
\[
\Delta_2:=D(f)\setminus\Delta_1.
\]
If $\Delta_1=\emptyset$ then $D(f)$ solely consists of points $\alpha^\star\in (2\N_0^n) \text{ with }f_{\alpha^\star}>0$ and the conditions \eqref{eq:weigh2}-\eqref{eq:weigh3} are thus satisfied for any choice of weights $\omega(\alpha^\star)>0$, $\alpha^\star\in D(f)$ fulfilling \eqref{eq:weigh1}. In the following let us thus consider only the case $\Delta_1\not=\emptyset$. If $\Delta_1=D(f)$ then due to $D(f)=\Delta_1\ \dot\cup\ \Delta_2$ one obtains $\Delta_2=\emptyset$ and one can define weights
\begin{equation}\label{eq:def_omega1}
\omega(\alpha^\star):=\frac{|f_{\alpha^\star}|}{\Theta(f,\lambda_f,\alpha^\star)}+\frac{\epsilon}{|\Delta_1|}\qquad\text{ for all }\alpha^\star\in\Delta_1.
\end{equation}
yielding

\[
\sum_{\alpha^\star\in D(f)}\omega(\alpha^\star)=\sum_{\alpha^\star\in\Delta_1}\omega(\alpha^\star)
\]
\[
=\sum_{\alpha^\star\in\Delta_1}\left( \frac{|f_{\alpha^\star}|}{\Theta(f,\lambda_f,\alpha^\star)}+\frac{\epsilon}{|\Delta_1|} \right)=\sum_{\alpha^\star\in\Delta_1} \frac{|f_{\alpha^\star}|}{\Theta(f,\lambda_f,\alpha^\star)}+\epsilon
\]
\[
=\sum_{\alpha^\star\in D(f)\,\cap\, (2\N_0^n)^c}\frac{|f_{\alpha^\star}|}{\Theta(f,\lambda_f,\alpha^\star)} +\sum_{\substack{\alpha^\star\in D(f)\, \cap\, 2\N_0^n \\ f_{\alpha^\star} <0}}\frac{|f_{\alpha^\star}|}{\Theta(f,\lambda_f,\alpha^\star)}\quad +\quad\epsilon
\]
\[
=\sum_{\alpha^\star\in D(f)\,\cap\, (2\N_0^n)^c}\frac{|f_{\alpha^\star}|}{\Theta(f,\lambda_f,\alpha^\star)} -\sum_{\alpha^\star\in D(f)\, \cap\, 2\N_0^n }\frac{\min\{0,f_{\alpha^\star}\}}{\Theta(f,\lambda_f,\alpha^\star)}\quad +\quad \epsilon\quad=\quad 1,
\]
where the second last equality is true due to $-\min\{0,f_{\alpha^\star}\}=|f_{\alpha^\star}|$ for all $\alpha^\star\in\Delta_1$ with $\alpha^\star\in 2\N_0^n$ and the last equality follows by \eqref{eq:epsilon}. The condition \eqref{eq:weigh1} follows. By \ref{eq:def_omega1} for each $\alpha^\star\in D(f)$ one has
\[
\omega(\alpha^\star)>\frac{|f_{\alpha^\star}|}{\Theta(f,\lambda_f,\alpha^\star)},
\]
which directly implies the condition \eqref{eq:weigh3} and due to $|f_{\alpha^{\star}}|=-f_{\alpha{^\star}}$ for $\alpha^{\star}\in 2\N_0^n$, also the condition \eqref{eq:weigh2} follows.

It remains to check the last case where $\Delta_1\not=\emptyset$ and $\Delta_2\not=\emptyset$. Here one can define 

\[
\omega(\alpha^\star):=\frac{|f_{\alpha^\star}|}{\Theta(f,\lambda_f,\alpha^\star)}+\frac{\epsilon}{2|\Delta_1|}\qquad\text{ for each }\alpha^\star\in\Delta_1.
\]
together with
\[
\omega(\alpha^\star):=\frac{\epsilon}{2|\Delta_2|}\qquad\text{ for each }\alpha^\star\in\Delta_2.
\]
By $D(f)=\Delta_1\ \dot\cup\ \Delta_2$ this yields
\[
\sum_{\alpha^\star\in D(f)}\omega(\alpha^\star)=\sum_{\alpha^\star\in\Delta_1}\omega(\alpha^\star)+\sum_{\alpha^\star\in\Delta_2}\omega(\alpha^\star)
\]
\[
=\sum_{\alpha^\star\in\Delta_1}\left( \frac{|f_{\alpha^\star}|}{\Theta(f,\lambda_f,\alpha^\star)}+\frac{\epsilon}{2|\Delta_1|} \right)+ \sum_{\alpha^\star\in\Delta_2} \frac{\epsilon}{2|\Delta_2|}
\]
\[
=\sum_{\alpha^\star\in\Delta_1} \frac{|f_{\alpha^\star}|}{\Theta(f,\lambda_f,\alpha^\star)}+\sum_{\alpha^\star\in\Delta_1} \frac{\epsilon}{2|\Delta_1|}+\sum_{\alpha^\star\in\Delta_2} \frac{\epsilon}{2|\Delta_2|} 
\]
\[
=\sum_{\alpha^\star\in\Delta_1} \frac{|f_{\alpha^\star}|}{\Theta(f,\lambda_f,\alpha^\star)}+\frac{\epsilon}{2}+\frac{\epsilon}{2}
\]
\[
=\sum_{\alpha^\star\in D(f)\,\cap\, (2\N_0^n)^c}\frac{|f_{\alpha^\star}|}{\Theta(f,\lambda_f,\alpha^\star)} +\sum_{\substack{\alpha^\star\in D(f)\, \cap\, 2\N_0^n \\ f_{\alpha^\star} <0}}\frac{|f_{\alpha^\star}|}{\Theta(f,\lambda_f,\alpha^\star)}\quad +\quad\epsilon
\]
\[
=\sum_{\alpha^\star\in D(f)\,\cap\, (2\N_0^n)^c}\frac{|f_{\alpha^\star}|}{\Theta(f,\lambda_f,\alpha^\star)} -\sum_{\alpha^\star\in D(f)\, \cap\, 2\N_0^n }\frac{\min\{0,f_{\alpha^\star}\}}{\Theta(f,\lambda_f,\alpha^\star)}\quad +\quad \epsilon\quad=\quad 1
\]
where the second last equality is true due to $-\min\{0,f_{\alpha^\star}\}=|f_{\alpha^\star}|$ for all $\alpha^\star\in\Delta_1$ with $\alpha^\star\in 2\N_0^n$ and due to $\min\{0,f_{\alpha^\star}\}=0$ for all $\alpha^\star\in\Delta_2$. The last equality follows by \eqref{eq:epsilon}. The property \eqref{eq:weigh1} is thus satisfied, and finally, it is easy to see that the conditions \eqref{eq:weigh2} and \eqref{eq:weigh3} are fulfilled as well analogous to the case $\Delta_2=\emptyset$ from above.
\qed

Now, by using Lemma \ref{lem:polyhed_suff}, we may restate Theorem \ref{the:suffdeg} as follows.

\bthe\label{the:suffdeg_restated}Let $f\in\R[x]$ be a polynomial satisfying the conditions \eqref{eq:c1}-\eqref{eq:c3} and let $\lambda_f$ be a map of minimal barycentric coordinates of $f$. If
\begin{equation}\label{eq:the_suffdeg_restated}
\sum_{\alpha^\star\in D(f)\,\cap\, (2\N_0^n)^c}\frac{|f_{\alpha^\star}|}{\Theta(f,\lambda_f,\alpha^\star)} -\sum_{\alpha^\star\in D(f)\, \cap\, 2\N_0^n }\frac{\min\{0,f_{\alpha^\star}\}}{\Theta(f,\lambda_f,\alpha^\star)}<1
\end{equation}
then $f$ is coercive on $\R^n$.
\ethe

Note that by applying the definition of the circuit number we can rewrite the left-hand side of the expression \eqref{eq:the_suffdeg_restated} as follows
\[
\sum_{\alpha^\star\in D(f)\,\cap\, (2\N_0^n)^c}\frac{|f_{\alpha^\star}|}{\Theta(f,\lambda_f,\alpha^\star)} -\sum_{\alpha^\star\in D(f)\, \cap\, 2\N_0^n }\frac{\min\{0,f_{\alpha^\star}\}}{\Theta(f,\lambda_f,\alpha^\star)}
\]
\[
=\sum_{\alpha^\star\in D(f)\,\cap\, (2\N_0^n)^c}  |f_{\alpha^\star}| \prod_{\alpha\in W_{\alpha^\star}(\lambda_f)}f_\alpha^{-\lambda_f(\alpha^\star,\alpha)} \prod_{\alpha\in W_{\alpha^\star}(\lambda_f)} \lambda_f(\alpha^\star,\alpha)^{\lambda_f(\alpha^\star,\alpha)}
\]

\[
-\sum_{\alpha^\star\in D(f)\, \cap\, 2\N_0^n }\min\{0,f_{\alpha^\star}\}\prod_{\alpha\in W_{\alpha^\star}(\lambda_f)}f_\alpha^{-\lambda_f(\alpha^\star,\alpha)} \prod_{\alpha\in W_{\alpha^\star}(\lambda_f)} \lambda_f(\alpha^\star,\alpha)^{\lambda_f(\alpha^\star,\alpha)} 
\]

which reveals that \eqref{eq:the_suffdeg_restated} can be viewed as strict posynomial inequality in variables $f_\alpha,\ \alpha\in V(f)$; $ -\min\{0,f_{\alpha^\star}\},\  \alpha^\star\in D(f)\,\cap\, (2\N_0^n)^c$ and $|f_{\alpha^\star}|,\ \alpha^\star\in D(f)\, \cap\, 2\N_0^n$. For more details on theory of posynomials and geometric programming we refer to \cite{Boyd}.

\bex\label{ex:compare2} Consider the homogenous bivariate quartic $f(x,y)=ax^4+bx^3y+cy^4$ with some parameter values $a,b,c\in\R\setminus\{0\}$. If $a,c\in\R_{>0}$, then $f$ fulfills the conditions \eqref{eq:c1}-\eqref{eq:c3}. If additionally to $a,c\in\R_{>0}$ the condition $b=0$ is fulfilled, then $f$ is gem regular and according to Theorem \ref{the:char} $f$ is coercive on $\R^2$. If $b\not=0$, then $f$ is gem degenerate with $D(f)=\{(3,1)\}$ and we obtain a unique map of minimal barycentric coordinates $\lambda_f$ of $f$ with
\[
\lambda_f((3,1),(4,0))=\frac{3}{4}\quad \lambda_f((3,1),(0,4))=\frac{1}{4}.
\]
which yields the circuit number
\[
\Theta(f,\lambda_f,(3,1))=\prod_{\alpha\in \{(4,0),(0,4)\}}\left(\frac{f_\alpha}{\lambda_f((3,1),\alpha)}\right)^{\lambda_f((3,1),\alpha)}=\left(\frac{a}{\frac{3}{4}}\right)^{\frac{3}{4}}\cdot \left(\frac{c}{\frac{1}{4}}\right)^{\frac{1}{4}}
\]
By Theorem \ref{the:suffdeg_restated} we obtain coercivity of $f$ if $a,c\in\R_{>0}$ and $b$ satisfy \eqref{eq:the_suffdeg_restated}, that is if
\begin{equation}
(a,b,c)\in\mathcal{C}_1:=\{(a,b,c)\in\R^3|\ \frac{|b|}{4\cdot 3^{-\frac{3}{4}}\cdot a^{\frac{3}{4}}\cdot c^{\frac{1}{4}}}<1\text{ and } a,c\in\R_{>0}\}\label{eq:coerc_3d}
\end{equation}
By Theorem \ref{thm:main} we obtain coercivity of $f$ if $a,c\in\R_{>0}$ and $b$ satisfy \eqref{eq:main}, that is if
\[
(a,b,c)\in\{(a,b,c)\in\R^3|\ a>\frac{3}{4}\cdot |b|,\  c\ >\ \frac{1}{4}\cdot |b|\text{ and } a,c\in\R_{>0}\}
\]
which is equivalent to 
\begin{equation}
(a,b,c)\in\mathcal{C}_2:=\{(a,b,c)\in\R^3|\ |b|<\min \,\{ \frac{4}{3}\cdot a, 4\cdot c\}\text{ and } a,c\in\R_{>0}\}.\label{eq:coerc_ed_polyhedral}
\end{equation}
\eex

\begin{figure}[htb]
\begin{subfigure}{.48\textwidth}
\centering
\includegraphics[width=1.0\linewidth]{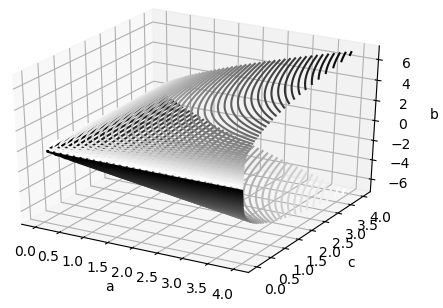}
\caption{Boundary of the convex set of coefficients $\mathcal{C}_1$ implying coercivity of $f$.}
\label{subfig:11}
\end{subfigure}
\hfill
\begin{subfigure}{.48\textwidth}
\centering
\includegraphics[width=1.0\linewidth]{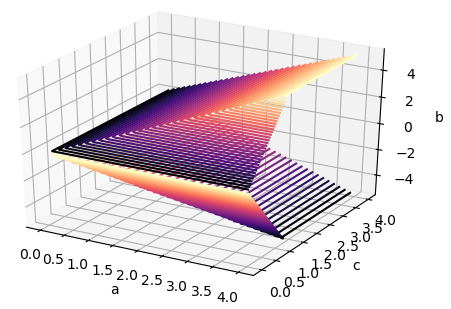}
\caption{Boundary of the convex set of coefficients $\mathcal{C}_2$ implying coercivity of $f$.}
\label{subfig:12}
\end{subfigure}
\hfill
\begin{subfigure}{.48\textwidth}
\centering
\includegraphics[width=1.0\linewidth]{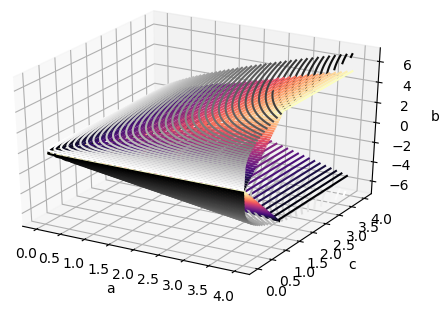}
\caption{Relationship between the sets $\mathcal{C}_1$ and $\mathcal{C}_2$.}
\label{subfig:21}
\end{subfigure}
\begin{subfigure}{.48\textwidth}
\centering
\includegraphics[width=1.0\linewidth]{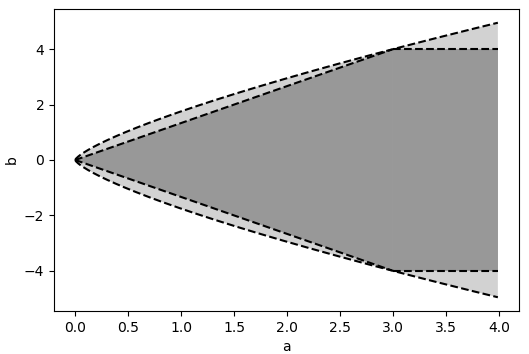}
\caption{Slices of sets $\mathcal{C}_1$ and $\mathcal{C}_2$ for parameter value $c=1$.}
\label{subfig:22}
\end{subfigure}
\caption{Subsets of polynomial coefficients $\mathcal{C}_1$ and $\mathcal{C}_2$ from Example \ref{ex:compare2}.}
\label{fig:ALL}
\end{figure}

In the latter Example \ref{ex:compare2} it becomes aparent that the subset of polynomial coefficients $\mathcal{C}_2\subset\R^3$ implying coercivity of $f$ by Theorem \ref{thm:main} is contained in the set $\mathcal{C}_1$ identified by Theorem \ref{the:suffdeg_restated}. A natural question arising in this context is whether this is always true, or whether there exist polynomials for which Theorem \ref{thm:main} identifies coefficients implying their coercivity which are not captured by Theorem \ref{the:suffdeg_restated}. Next Example \ref{ex:compare} shows that the latter is true in general. Moreover we will show that the inclusion $\mathcal{C}_2\subseteq\mathcal{C}_1$ in Example \ref{ex:compare2} is a consequence of the following result, which sharpens the existing result on necessary conditions for coercivity for a class of gem irregular polynomials (see Theorem 2.29 and Remark 3.6 in \cite{BS}).

\bthe[Coercivity characterization for circuit polynomials]\label{the:Char_Circ} Let $f\in\R[x]$ of the form
\[
f(x)=\sum_{j=0}^rf_{\alpha(j)}x^{\alpha(j)}+f_{\alpha^\star}x^{\alpha^\star}
\] 
with $r\leq n$ be a circuit polynomial according to Definition \ref{defi:circuit_polynomial} and let further $\lambda_f$ denote the unique map of minimal barycentric coordinates of $f$. Then $f$ is coercive on $\R^n$ if and only if one of the following conditions is fulfilled:
\begin{itemize}
\item[(a)] $r=n$ and $f$ fulfills conditions \eqref{eq:c1}--\eqref{eq:c3}
\item[(b)] $r=n-1$ and $f$ fulfills conditions \eqref{eq:c1}--\eqref{eq:c3} together with
\begin{equation}\label{eq:strict1}
f_{\alpha^\star}\ >\ -\Theta(f,\lambda_f,\alpha^\star)\ \text{ if }\ \alpha^\star\in 2\N_0^n
\end{equation}
and 
\begin{equation}\label{eq:strict2}
|f_{\alpha^\star}|\ <\ \Theta(f,\lambda_f,\alpha^\star)\ \text{ else.}
\end{equation}
\end{itemize}
\ethe

\proof \textbf{"$\Rightarrow$"} For the first direction let $f$ be coercive on $\R^n$. By Theorem \ref{the:necessary_general} polynomial $f$ fulfills the conditions \eqref{eq:c1}--\eqref{eq:c3}. Due to \eqref{eq:c3} the inequality $|V(f)|\geq n$ holds and since for the set of vertices at infinity $V(f)$ of $f$ the inclusion $V(f)\subseteq\text{Vert}(\new(f))$ holds, we obtain $|\text{Vert}(\new(f))|\geq n$. By Definition \ref{defi:circuit_polynomial} of a circuit polynomial we have that the set $\text{Vert}(\new(f))=\{\alpha(0),\dots,\alpha(r)\}\subseteq\R^n$ is affinely independent and thus $|\text{Vert}(\new(f))|\leq n+1$ holds necessarily. Combining both inequalities yields the bounds $n\leq |\text{Vert}(\new(f))| \leq n+1$ and thus $r\in\{n-1,n\}$ follows. Since the conditions \eqref{eq:c1}-\eqref{eq:c3} are fulfilled independently of the exact value $r\in\{n-1,n\}$, the assertion $(a)$ is already shown. For concluding the proof of the first direction it remains thus to show that additionally to conditions \eqref{eq:c1}-\eqref{eq:c3} also the strict inequalities \eqref{eq:strict1} and \eqref{eq:strict2} hold in case $r=n-1$. In fact, due to condition \eqref{eq:c3} all vectors in the set $\text{Vert}(\new(f))=\{\alpha(0),\dots,\alpha(n-1)\}$ can be written in the form $\alpha(j)=2k_{i(j)}e_{i(j)}$ with some $k_{i(j)}\in\N,\ j=0,\dots,n-1$ and with standard unit vectors $e_{i(j)}$ where 
\[
i:\{0,\dots,n-1\}\to\{1,\dots,n\}
\] 
denotes the corresponding bijective index map. The latter together with the representation of the vector $\alpha^\star$ from the Definition \ref{defi:circuit_polynomial} of circuit polynomial $f$ yield $\alpha^\star\in V^c(f)\,\cap\, G$ with a simplicial face $G:=\conv\{\alpha(0),\dots,\alpha(n-1)\}$ of $\new_\infty(f)$ fulfilling $0\notin G$. Thus gem irregularity of $f$ with $D(f)=\{\alpha^\star\}$ follows and application of Theorem 2.29 and Remark 3.6 from \cite{BS} yields that coercivity of $f$ implies
\begin{equation}\label{eq:non_strict1}
f_{\alpha^\star}\ \geq\ -\Theta(f,\lambda_f,\alpha^\star)\ \text{ if }\ \alpha^\star\in 2\N_0^n
\end{equation}
and 
\begin{equation}\label{eq:non_strict2}
|f_{\alpha^\star}|\ \leq\ \Theta(f,\lambda_f,\alpha^\star)\ \text{ else.}
\end{equation}

In order to conclude the proof of the first direction, it hence suffices to show that in both cases \eqref{eq:non_strict1} and \eqref{eq:non_strict2} the respective equality cannot occur under the presence of coercivity of $f$. For this aim let us first consider the case $\alpha^\star\in 2\N_0^n$ and assume that \eqref{eq:non_strict1} is fulfilled with equality, that is, $f_{\alpha^\star}=-\Theta(f,\lambda_f,\alpha^\star)$. Then the circuit polynomial $f$ can be written in the form
\[
f(x)\quad =\quad \sum_{j=0}^{n-1}f_{\alpha(j)}x^{\alpha(j)}\quad-\quad\Theta(f,\lambda_f,\alpha^\star)\cdot x^{\alpha^\star}
\]
\[
=\quad \sum_{j=0}^{n-1}f_{\alpha(j)}x^{\alpha(j)}\quad -\quad \prod_{\alpha\in W_{\alpha^\star}(\lambda_f)}\left(\frac{f_\alpha}{\lambda_f(\alpha^\star,\alpha)}\right)^{\lambda_f(\alpha^\star,\alpha)} \cdot x^{\alpha^\star}
\]
\[
=\quad \sum_{j=0}^{n-1}f_{\alpha(j)}x^{\alpha(j)}\quad -\quad \prod_{j=0}^{n-1}\left(\frac{f_{\alpha(j)}}{\lambda_j}\right)^{\lambda_j} \cdot x^{\sum_{j=0}^{n-1}\lambda_j\alpha(j)}
\]
\begin{equation}\label{eq:AGM}
=\quad \sum_{j=0}^{n-1}\lambda_j\left(\frac{f_{\alpha(j)}x^{\alpha(j)}}{\lambda_j}\right) \quad-\quad \prod_{j=0}^{n-1}\left(\frac{f_{\alpha(j)}}{\lambda_j}\right)^{\lambda_j} \cdot x^{\sum_{j=0}^{n-1}\lambda_j\alpha(j)}
\end{equation}
where for the second equality we just apply the definition of the circuit number and for the third equality we use the property iv) from Definition \ref{defi:circuit_polynomial} of the circuit polynomial together with the property $W_{\alpha^\star}(\lambda_f)=\{\alpha(0),\dots,\alpha(n-1)\}$ of the minimal vertex representation of $\alpha^\star$ corresponding to $\lambda_f$. Due to the conditions \eqref{eq:c1} and \eqref{eq:c2} as well as property iv) from Def. \ref{defi:circuit_polynomial}, we further  have
\[
\frac{f_{\alpha(j)}x^{\alpha(j)}}{\lambda_j}\geq 0\quad\text{for all }x\in\R^n,\ j=0,\dots,n-1
\] 
and a direct application of the weighted arithmetic-geometric-mean inequality on the expression \eqref{eq:AGM} yields $f(x)\geq 0$ for all $x\in\R^n$ as well as the property $f(x)=0$ with $x\not= 0$ if and only if
\[
\frac{f_{\alpha(j)}x^{\alpha(j)}}{\lambda_j}=t\quad \text{for all $j=0,\dots,n-1$}
\]
with some fixed constant $t>0$. The property $f(x(t))\equiv 0$ is fulfilled, for instance, along the one-dimensional manifold $x(t):\R_{>0}\to(\R_>)^n$ of the form
\begin{align}\label{eq:manifold}
    x(t) &=\begin{bmatrix}
           \left(t\cdot\frac{\lambda_{i^{-1}(1)}}{f_{\alpha(i^{-1}(1))}}\right)^{\frac{1}{2k_{i^{-1}(1)}}} \\
           \vdots \\
           \left(t\cdot\frac{\lambda_{i^{-1}(n)}}{f_{\alpha(i^{-1}(n))}}\right)^{\frac{1}{2k_{i^{-1}(n)}}}
         \end{bmatrix}.
  \end{align}
Since $f_{\alpha(i^{-1}(j))}>0$ is fulfilled for each $j=1,\dots,n$ due to condition \eqref{eq:c2} and furthermore $\lambda_{i^{-1}(j)}>0$ and $k_{i^{-1}(j)}\in\N$ for each $j=1,\dots,n$ we obtain $\left(t\cdot\frac{\lambda_{i^{-1}(j)}}{f_{\alpha(i^{-1}(j))}}\right)^{\frac{1}{2k_{i^{-1}(j)}}}\to +\infty$ for $t\to +\infty$ and for each $j=1,\dots,n$. This yields $\|x(t)\|\to +\infty$ for $t\to +\infty$ with $f(x(t))\equiv 0$ contradicting the coercivity of $f$ on $\R^n$ and the strict inequality \eqref{eq:strict1} follows. 

For the other case $\alpha^\star\in (2\N_0^n)^c$ we have to show that the equality $f_{\alpha^\star}=|\Theta(f,\lambda_f,\alpha^\star)|$ can not occur if $f$ is coercive on $\R^n$. Here again, if $f_{\alpha^\star}=-\Theta(f,\lambda_f,\alpha^\star)$, the same line of argumentation can be applied as above. On the other hand, if $f_{\alpha^\star}=\Theta(f,\lambda_f,\alpha^\star)$, analogous to \eqref{eq:AGM} we obtain
\begin{equation}\label{eq:AGM2}
f(x)=\quad \sum_{j=0}^{n-1}\lambda_j\left(\frac{f_{\alpha(j)}x^{\alpha(j)}}{\lambda_j}\right) \quad+\quad \prod_{j=0}^{n-1}\left(\frac{f_{\alpha(j)}}{\lambda_j}\right)^{\lambda_j} \cdot x^{\sum_{j=0}^{n-1}\lambda_j\alpha(j)},
\end{equation}
and since $\alpha^\star\in (2\N_0^n)^c$ there exists some $s\in\{1,\dots,n\}$ with an odd entry $\alpha^\star_s\in (2\N_0)^c$. Using the map $x(\cdot)$ from \eqref{eq:manifold} let us define $\tilde{x}(t):\R_{>0}\to\R^n$ by 
\[\tilde{x}_j(t):= \begin{cases} 
      x_j(t) & \text{if }j\in\{1,\dots,n\}\setminus\{s\} \\
      -x_j(t) & j=s 
   \end{cases}
\]
which fulfills the property $\|\tilde{x}(t)\|\to +\infty$ for $t\to +\infty$ analogous to $x(t)$. Inserting $\tilde{x}(t)$ into \eqref{eq:AGM2} yields
\[
f(\tilde{x}(t))=\quad \sum_{j=0}^{n-1}\lambda_j\left(\frac{f_{\alpha(j)}\left(\tilde{x}(t)\right)^{\alpha(j)}}{\lambda_j}\right) \quad+\quad \prod_{j=0}^{n-1}\left(\frac{f_{\alpha(j)}}{\lambda_j}\right)^{\lambda_j} \cdot \left(\tilde{x}(t)\right)^{\sum_{j=0}^{n-1}\lambda_j\alpha(j)}
\]
\[
=\quad \sum_{j=0}^{n-1}\lambda_j\left(\frac{f_{\alpha(j)}\left(x(t)\right)^{\alpha(j)}}{\lambda_j}\right) \quad-\quad \prod_{j=0}^{n-1}\left(\frac{f_{\alpha(j)}}{\lambda_j}\right)^{\lambda_j} \cdot \left(x(t)\right)^{\sum_{j=0}^{n-1}\lambda_j\alpha(j)}\equiv 0
\]
for all $t>0$, where the second equality holds due to eveness property $\alpha(j)\in 2\N_0^n$ for each $j=0,\dots,n-1$ as well as due to the property $(\tilde{x}_s(t))^{\alpha^\star_s}=-(x_s(t))^{\alpha^\star_s}$ and the last equality follows from the weighted arithmetic-geometric mean inequality. For the case $\alpha^\star\in (2\N_0^n)^c$ the property $f_{\alpha^\star}=|\Theta(f,\lambda_f,\alpha^\star)|$ thus always enables a construction of a one-dimensional manifold $\tilde{x}(t)$ with $\|\tilde{x}(t)\|\to +\infty$ for $t\to +\infty$ satisfying $f(\tilde{x}(t))\equiv 0$. This contradicts the coercivity of $f$ on $\R^n$ and the strict inequality \eqref{eq:strict2} follows.

\textbf{"$\Leftarrow$"} For the other direction it suffices to show that in case $r=n$ the circuit polynomial $f$ satisfying \eqref{eq:c3} is gem regular and that in case $r=n-1$ the circuit polynomial $f$ satisfying \eqref{eq:c3} is gem irregular with $D(f)=\{\alpha^\star\}$. Then namely a direct application of Theorem \ref{the:char} for case (a) and of Theorem \ref{the:suffdeg_restated} for case (b) yields coercivity of $f$. 
Since in the proof of the first direction "$\Rightarrow$", we have already shown that in case $r=n-1$ the property \eqref{eq:c3} implies gem-irregularity of the circuit polynomial $f$ with $D(f)=\{\alpha^\star\}$, it only remains to tackle the case $r=n$. For this aim let the circuit polynomial $f$ with $r=n$ be given. We want to show that under the presence of \eqref{eq:c3} $f$ is gem regular. Since $\{\alpha(0),\dots,\alpha(n)\}$ are affinely independent, Newton polytope of $f$ 
\[
\new(f)=\conv \{\alpha(0),\dots,\alpha(n)\}
\] 
is a full-dimensional polytope containing the exponent vector $\alpha^\star$ as its inner point. But then also the Newton polytope at infinity of $f$
\begin{equation}\label{eq:Newton_vertices}
\new_\infty(f)=\conv \{0,\alpha(0),\dots,\alpha(n)\}
\end{equation}
fulfilling $\new(f)\subseteq\new_{\infty}(f)$ is a full-dimensional polytope which contains the exponent vector $\alpha^\star$ as its inner point. Due to full-dimensionality of the Newton polytope at infinity $\new_\infty(f)$ we obtain for its set of vertices the property $|\text{Vert}(\new_\infty(f))|\geq n+1$. Due to \eqref{eq:Newton_vertices} the upper bound $|\text{Vert}(\new_\infty(f))|\leq n+2$ is always fulfilled and $|\text{Vert}(\new_\infty(f))|\in\{n+1,n+2\}$ follows. In case $|\text{Vert}(\new_\infty(f))|=n+2$ one has $V(f)=\{\alpha(0),\dots,\alpha(n)\}$ and since $V^c(f)=A(f)\setminus V(f)=\{\alpha^\star\}$, the only possible candidate for a gem degenerate exponent vector of $f$ is the exponent vector $\alpha^\star$. Since $\alpha^\star$ is an inner point of $\new_\infty(f)$, it can not be contained in any proper face of $\new_\infty(f)$ which does not include the origin. Thus, by Definition \ref{def:Fregular}, the exponent vector $\alpha^\star$ is not gem degenerate and gem regularity of $f$ follows. In case $|\text{Vert}(\new_\infty(f))|=n+1$, due to $0\in\text{Vert}(\new_\infty(f))$ and condition \eqref{eq:c3}, precisely $n$ points from the set $\{\alpha(0),\dots,\alpha(n)\}$ are vertices of $\new_\infty(f)$ and they are of the form $2k_ie_i$ with some $k_i\in\N$ for $i=1,\dots,n$ (w.l.o.g. assume these points are $\alpha(1),\dots,\alpha(n)$). Exponent vector $\alpha(0)$ thus eihter fulfills $\alpha(0)=0$ or $\alpha(0)$ is not a vertex of $\new_\infty(f)$. In case $\alpha(0)=0$ the only possible candidade for a gem degenerate exponent vector of $f$ is again the exponent vector $\alpha^\star$. Using the same line of argumentation for $\alpha^\star$ as above, gem-regularity of $f$ follows. In case $\alpha(0)$ is not a vertex of $\new_\infty(f)$, we have two candidates for gem degenerate exponent vectors of $f$ - the vectors $\alpha(0)$ and $\alpha^\star$. As we have already seen, $\alpha^\star$ as an inner point of $\new_\infty(f)$ can't be gem degenerate exponent vector of $f$, so it only remains to analyze $\alpha(0)$ for gem degeneracy. We have $\new_\infty(f)=\conv \{0,2k_1e_1,\dots,2k_ne_n\}$ and the only possibility for $\alpha(0)$ to be contained in some non-trivial face of $\new_\infty(f)$ which does not include the origin is given only if $\alpha(0)\in\conv\{\alpha(1),\dots,\alpha(n)\}$. This contradicts the assumption of affine independence of vectors $\{\alpha(0),\dots,\alpha(n)\}$ and the gem regularity of $f$ follows.
\qed

\brem In Example \ref{ex:compare2} the polynomial $f$ is a circuit polynomial with $r=n-1$. According to Theorem \ref{the:Char_Circ} the coercivity of $f$ is characterized by conditions \eqref{eq:c1}--\eqref{eq:c3} and the strict inequality \eqref{eq:strict2}. Theorem \ref{the:Char_Circ} implies that the subset of coefficients $\mathcal{C}_1$ as identified by Theorem \ref{the:suffdeg_restated} and depicted in Figure \ref{fig:ALL} is complete in the sense that it contains all polynomial coefficients for which $f$ is coercive on $\R^2$. Even if in Example \ref{ex:compare2} our Theorem \ref{thm:main} identified only a strict (polyhedral) subset $\mathcal{C}_2$ of the set $\mathcal{C}_1$, in the next example we will see that in general Theorem \ref{thm:main} can identify coercive polynomials which are not captured by Theorem \ref{the:suffdeg_restated}. 
\erem

\bex\label{ex:compare} Consider the homogenous bivariate quartic $g(x,y)=x^4+ax^3y+bxy^3+y^4$ with some parameter values $a,b\in\R$. One obtains $V(g)=\{(4,0),(0,4)\}$ and $D(g)=\{(3,1),(1,3)\}$ which leads to unique map of minimal barycentric coordinates $\lambda_g$ of $g$ with $\lambda_g((1,3),(0,4))=\lambda_g((3,1),(4,0))=\frac{3}{4}$ and $\lambda_g((1,3),(4,0))=\lambda_g((3,1),(0,4))=\frac{1}{4}$. Computing the circuit numbers corresponding to all gem degenerate exponent vectors of $g$ yields
\[
\Theta(g,\lambda_g,(3,1))=\prod_{\alpha\in \{(4,0),(0,4)\}}\left(\frac{g_\alpha}{\lambda_g((3,1),\alpha)}\right)^{\lambda_g((3,1),\alpha)}=\left(\frac{1}{\frac{3}{4}}\right)^{\frac{3}{4}}\cdot \left(\frac{1}{\frac{1}{4}}\right)^{\frac{1}{4}}=4\cdot 3^{-\frac{3}{4}}
\]
\[
\Theta(g,\lambda_g,(1,3))=\prod_{\alpha\in \{(4,0),(0,4)\}}\left(\frac{g_\alpha}{\lambda_g((1,3),\alpha)}\right)^{\lambda_g((1,3),\alpha)}=\left(\frac{1}{\frac{1}{4}}\right)^{\frac{1}{4}}\cdot \left(\frac{1}{\frac{3}{4}}\right)^{\frac{3}{4}}=4\cdot 3^{-\frac{3}{4}}.
\]

According to Theorem \ref{the:suffdeg_restated} one obtains that  the polynomial $g$ is coercive on $\R^2$ for parameter values $a,b\in\R$ satisfying 
\begin{equation}\label{eq:pol_g_coerc_old}
\frac{|g_{(3,1)}|}{\Theta(g,\lambda_g,(3,1))}+\frac{|g_{(1,3)}|}{\Theta(g,\lambda_g,(1,3))}=\frac{|a|}{4\cdot 3^{-\frac{3}{4}}}+\frac{|b|}{4\cdot 3^{-\frac{3}{4}}}<1.
\end{equation}
On the other hand, according to Theorem \ref{thm:main}, polynomial $g$ is coercive on $\R^2$ for parameter values $a,b\in\R$ satisfying the inequalities
\begin{equation}\label{eq:pol_g_coerc_new}
1>\frac{3}{4}\cdot|a|+\frac{1}{4}\cdot|b|\quad\text{and}\quad 1>\frac{1}{4}\cdot|a|+\frac{3}{4}\cdot|b|.
\end{equation}

\begin{figure}[htb]
\begin{subfigure}{.5\linewidth}
\centering
\includegraphics[width=1.0\linewidth]{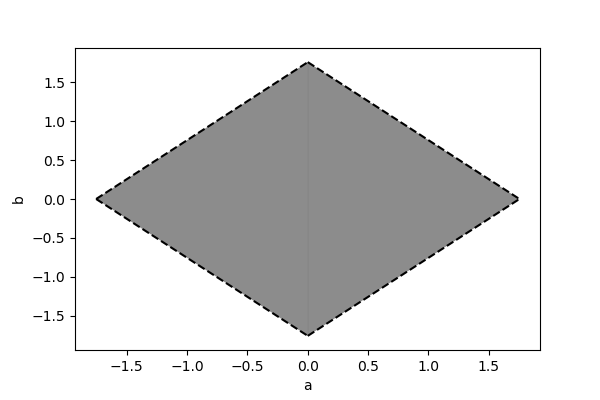}
\caption{The rhombic area depicts coefficients $(a,b)\in\R^2$ satisfying \eqref{eq:pol_g_coerc_old}}
\label{fig:sub1}
\end{subfigure}%
\begin{subfigure}{.5\linewidth}
\centering
\includegraphics[width=1.0\linewidth]{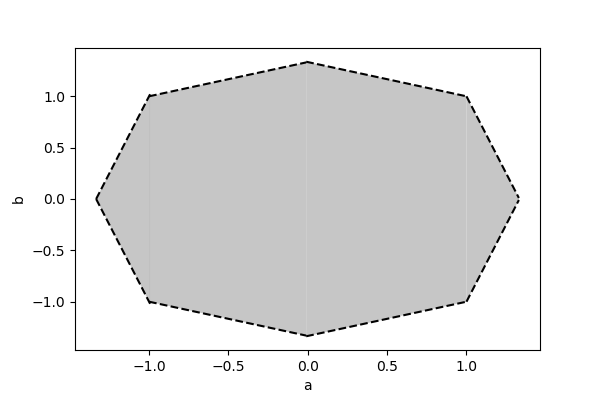}
\caption{The hexagonal area depicts coefficients $(a,b)\in\R^2$ satisfying \eqref{eq:pol_g_coerc_new}}
\label{fig:sub2}
\end{subfigure}\\[1ex]
\begin{subfigure}{\linewidth}
\centering
\includegraphics[width=0.5\linewidth]{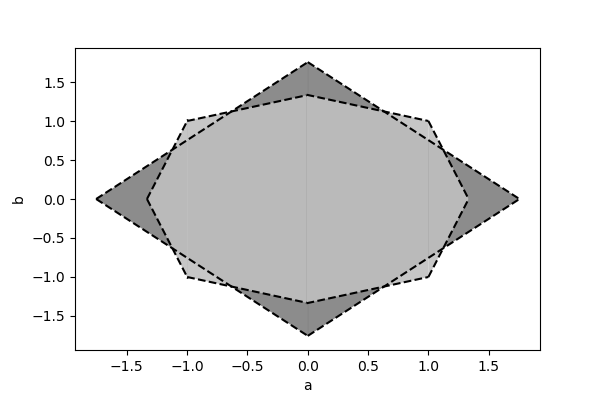}
\caption{Relation between the two - the rhombic and the hexagonal area}
\label{fig:sub3}
\end{subfigure}
\caption{Polyhedral subsets from Example \ref{ex:compare}.}
\label{fig:2}
\end{figure}

The four light-shaded triangular areas lying outside the intersection of the rhombic and the hexagonal area depicted in Figure \ref{fig:sub3} represent those coefficients $(a,b)\in\R^2$ of the polynomial $g$ for which $g$ is coercive according to Theorem \ref{thm:main} and for which Theorem  \ref{the:suffdeg_restated} does not imply coercivity of $g$. The four dark-shaded quadrilateral areas lying outside the intersection of the rhombic and the hexagonal area represent those coefficients $(a,b)\in\R^2$ for which the reverse is true.
\eex

%\begin{remark}
%Example \ref{ex:compare} shows that the sufficiency conditions for coercivity from Theorem \ref{the:suffdeg_restated} are in general independent from those given in Theorem \ref{thm:main}.
%\end{remark}

\section{Final remarks}\label{sec:fin}
With Theorem \ref{thm:main} we provide new conditions on Newton polytopes at infinity and on polynomial coefficients implying coercivity of general polynomials, which are independent from those identified in \cite{BS}. In fact, as shown in Example \ref{ex:compare}, our new conditions can in some cases guarantee coercivity for polynomials even if the other known sufficiency conditions are not satisfied and vice versa. With Characterization Theorem \ref{the:Char_Circ} we furthermore enlarge the class of polynomials beyond the class of gem regular polynomials introduced in \cite{BS} for which a characterization of their coercivity can be given using conditions involving Newton polytopes. As for the cone of Sum-of-Nonnegative-Circuits (SONC), the Characterization Theorem \ref{thm:IliWolff} forms a theoretical basis for developing algorithmically tractable nonnegativity certificates for polynomials, it could be interesting to consider the cone of Sum-of-Coercive-Circuits (SOCC) in light of the Characterization Theorem \ref{the:Char_Circ} accordingly. Although it was recently shown in \cite{AAA} that checking the coercivity property even for low degree polynomial instances is NP-hard, it would be still of practical interest to identify and describe a broader class of coercive polynomials (such as e.g. SOCC) for which the coercivity could be verified in some systematic and tractable way. In \cite{GhaMar,L} sufficiency conditions for polynomials for being sum of squares are identified which are linear in polynomial coefficients and are hence of alike nature as those identified in Theorem \ref{thm:main}. This could be used to further analyze the structural differences between the cone of coercive polynomials and the cones of sum of squares or non-negative polynomials. We leave these aspects for future research.

\section*{Acknowledgments} The authors are grateful to Thorsten Theobald for his support and for fruitful discussions on the subject of this article.

\end{document}